\input amstex
\input epsf
\input pstricks
\documentstyle{amsppt}
\NoBlackBoxes
\nologo
\def\ac{\text{\rm arcsinh}}

\def\eps{\varepsilon}

\def\grad{\nabla}
\def\H{\Cal H}
\def\Hs{\Cal H^{*}}

\def\L{\Lambda}
\def\M{\mu^{\wh\beta}}
\def\norm#1{|#1|}
\def\N{\Bbb N}

\def\R{\Bbb R}
\def\vphi{\varphi}
\def\s{\sigma}
\def\T{\tau_\beta}
\def\ths{\theta^{*}}
\def\us{{u^{*}}}
\def\xs{x^{*}}

\def\Z{\Bbb Z}

\def\Exp#1{\exp\left(#1\right)}

\def\da{\downarrow}

\def\lra{\leftrightarrow}
\def\ra{\rightarrow}
\def\ua{\uparrow}
\def\wh{\widehat}

\def\rond#1{\smash{#1%
\raise 7.9pt%
\hbox{%
$\kern -9.4pt{\scriptstyle \circ}$%
}}}


\def\ea#1{\eqalign{#1}}

\def\marginal#1{\strut\setbox0=%
      \vtop{\hsize=15mm
            \fiverm
            \textfont0=\scriptfont0 \textfont1=\scriptfont1
            \textfont=\scriptfont2 \textfont3=\scriptfont3
            \parindent=0pt\baselineskip=7pt
             \raggedright\rightskip=1em plus 2em\hfuzz=.5em\tolerance=9000
             \overfullrule=0pt
             #1}%
      \dimen0=\ht0
      \advance\dimen0 by \dp\strutbox
      \ht0=0pt\dp0=\dimen0
      \vadjust{\kern-\dimen0\moveleft\wd0\box0}%
      \ignorespaces}

\def\noi{\noindent}

\def\q{\quad}

\def\qq{\qquad}

\def\today{\noi\number\day
\space\ifcase\month\or
  January\or February\or March\or April\or May\or June\or
  July\or August\or September\or October\or November\or December\fi
  \space\number\year}



\newcount\spheadno
\global\spheadno=0
\def\nextspheadno{\global\advance\spheadno by1 %
\the\spheadno %
\global\headno=0%
\global\subheadno=0%
\global\subsubheadno=0%
\global\cno=0%
\global\formno=0}

\def\spheadlabel#1{\edef #1{\the\spheadno}}

\newcount\headno
\global\headno=0
\def\nextheadno{\global\advance\headno by1 %
\the\headno %
\global\subheadno=0%
\global\subsubheadno=0}
\def\headnum{\nextheadno.\ }
\def\headlabel#1{\edef #1{\the\headno}}

\newcount\subheadno
\global\subheadno=0
\def\nextsubheadno{\global\advance\subheadno by 1 %
\the\headno.%
\the\subheadno%
\global\subsubheadno=0}
\def\subheadnum{\nextsubheadno. }
\def\subheadlabel#1{\edef #1{\the\headno.\the\subheadno}}

\newcount\subsubheadno
\global\subsubheadno=0
\def\nextsubsubheadno{\global\advance\subsubheadno by 1 %
\the\headno.\the\subheadno.\the\subsubheadno}

\def\subsubheadlabel#1{\edef #1{\the\subsubheadno}}

\newcount\formno
\newcount\cno
\def\nextno{\global\advance\cno by1 \the\cno }
\newcount\figno
\def\nextfigno{\global\advance\figno by1 \the\figno }

\def\figlabel#1{\edef #1{\the\figno}}

\global\formno=0
\def\nextformno{\global\advance\formno by1 \the\formno }

\def\blankpage{\vfil\nopagenumbers\eject\line{}%
\vfil\eject\global\advance\pageno
by -1\footline={\hss\tenrm\folio\hss}}

\def\eqnum{\tag{\nextformno}}
\def\eqlabel#1{%
\edef%
#1{\the\formno}%
}
\def\eqref#1{($#1$)}
\def\procnum{\nextno}
\def\proclabel#1{
\edef #1{\the\cno}}
\def\procref#1{$#1$}
\newcount\refno
\global\refno=0
\def\nextrefno{\global\advance\refno by 1 }
\nextrefno\edef\ALEXANDER{\the\refno}
\nextrefno\edef\BARROUCHMcCOYTRACYWU{\the\refno}
\nextrefno\edef\McCOYTRACY{\the\refno}
\nextrefno\edef\CAMIANEWMAN{\the\refno}
\nextrefno\edef\CamiaNewmanII{\the\refno}
\nextrefno\edef\CAMPANINOIOFFEVELENIKI{\the\refno}
\nextrefno\edef\CAMPANINOIOFFEVELENIKII{\the\refno}
\nextrefno\edef\CAMPANINOIOFFEVELENIKIII{\the\refno}
\nextrefno\edef\CERFMESSIKH{\the\refno}
\nextrefno\edef\CHAYESCHAYESCAMPANINO{\the\refno}
\nextrefno\edef\McCoyWu{\the\refno}
\nextrefno\edef\DZ{\the\refno}
\nextrefno\edef\GRIM{\the\refno}
\nextrefno\edef\HSU{\the\refno}  
\nextrefno\edef\IOFFE{\the\refno}
\nextrefno\edef\KASTELEYN{\the\refno}
\nextrefno\edef\KESTEN{\the\refno}
\nextrefno\edef\KRAMERSWANNIER{\the\refno}
\nextrefno\edef\MESSIKH{\the\refno}
\nextrefno\edef\MONTROLL{\the\refno} 
\nextrefno\edef\PFISTER{\the\refno}
\nextrefno\edef\SCHRAMM{\the\refno}
\nextrefno\edef\SIMON{\the\refno}
\nextrefno\edef\SMIRNOV{\the\refno}
\nextrefno\edef\SMIRNOVWERNER{\the\refno}
\nextrefno\edef\TRACY{\the\refno}
\nextrefno\edef\WATSON{\the\refno}

\topmatter

\title
The surface tension near criticality of the 2d-Ising model 
\endtitle

\author
R. J. Messikh
\endauthor

\affil
Chair of Stochastic Modeling, EPFL
\endaffil
\address
EPFL SB IMA CMOS INR 032 Station 14, 1015 Lausanne, Switzerland.
\endaddress
\email
redajurg.messikh\@epfl.ch
\endemail

\thanks
We would like to thank Rapha\"el Cerf for proposing the problem and for many pleasant and fruitful discussions. 
\endthanks

\date
20 October 2006
\enddate

\abstract 
For the two dimensional Ising model, we construct the adequate surface tension near criticality. 
The latter quantity has been shown \cite{\CERFMESSIKH} to play a central role in the study of phase coexistence in a joint limit where the temperature approaches the critical point from below and simultaneously the size of the system increases fast enough.
\endabstract
\endtopmatter
\document

\head\headnum {Introduction}\endhead
The present paper is devoted to a fundamental quantity related to the phenomenon of phase coexistence, namely the {\sl surface tension}.
In particular, we are interested in its behavior when we approach the critical point.

\noi There is a general belief that the surface tension and other quantities become {\sl isotropic} near the critical point.
Indeed, a wide family of statistical physics models defined on discrete lattices and taken at the critical point are conjectured to converge to a conformally invariant measure when the lattice mesh size goes to zero. 
In the particular case of independent site percolation on the triangular lattice, this conjecture has been made rigorous in \cite{\SMIRNOVWERNER} and \cite{\CAMIANEWMAN}. 
There, the continuous conformally invariant measure is described by the Schramm-Loewner Evolution \cite{\SCHRAMM}. 

\noi Critical phenomena also appear when one is not strictly at the critical point. 
When considering the model "near criticality", that is when one takes the thermodynamical limit and simultaneously approaches the critical  point, then  the asymptotic behavior of relevant quantities  is also influenced by critical phenomena. 
Among these joint limits, there is a special one, sometimes called the "scaling limit", at the threshold between two different behaviors. 

\noi On one side, if the temperature approaches the critical point fast enough, the model behaves as if the temperature is exactly at the critical point. 
In the case of independent percolation, such a behavior has been proved by Kesten in \cite{\KESTEN}. 

\noi On the other side, when the temperature approaches the critical point slowly enough then the influence of criticality is different.  
In this case,  the behavior of the model can be described by non-critical phenomena that are altered by criticality. 
Indeed in \cite{\CERFMESSIKH}, it has been proved that the 2d-Ising model taken at sub-critical temperatures at the vicinity of the critical point and in boxes that are large enough, still exhibits phase coexistence. 
In such regimes, the Wulff crystal persists but criticality washes-out the anisotropy  inherited from the geometry of the lattice and reduces the Wulff shape to an ordinary circle.

\noi The paper \cite{\CERFMESSIKH} required the construction of an adequate joint-limit surface tension. 
The proof of this construction has been sketched in \cite{\MESSIKH}, using heavily the explicit computations techniques of \cite{\McCoyWu}. 
In the present paper, we provide a construction of the {\sl surface tension near criticality} without adding an extra layer of explicit computations. 
Our results are, for the time being, restricted to the two dimensional Ising model on the square lattice. 
Indeed, the starting point of our analysis is the beautiful exact formula for the {\sl fixed temperature surface tension}. 

Using the duality property of Kramers-Wannier \cite{\KRAMERSWANNIER}, the surface tension of the 2d-Ising model at an inverse temperature $\beta>\beta_c$ can be defined from the asymptotic behavior of the two point function at the dual inverse temperature $\wh\beta<\beta_c$ \cite{\PFISTER}. 
To be more precise, for each  $\wh\beta<\beta_c$, let us denote by $\mu^{\wh\beta}$ the unique infinite volume Ising measure on the spin configurations $\s\in\{-1,1\}^{\Z^2}$.  
Then for each $x\in\Z^2$ and at each  $\beta>\beta_c$, the surface tension is the function defined by the following limit 
$$
\T(x)=-\lim_{n\ra\infty}{1\over n}\log\mu^{\wh\beta}[\s(0)\s(nx)],\eqnum
$$\eqlabel{\eqdefsurften} 
where $\wh\beta<\beta_c$ is related to $\beta>\beta_c$ by the duality relation
$$
\sinh(2\beta)\sinh(2{\wh\beta})=1.
$$
It is well known \cite{\PFISTER} that the function $\T$ can be extended continuously into a norm on $\R^2$. 

\noi An important particularity of the 2d-Ising model is its relation to the dimer model \cite{\KASTELEYN}. 
Indeed, Kasteleyn discovered that the partition function of the Ising model can be represented as the generating function of a dimer model. This permitted to give an explicit formula for the partition function. 

\noi Later, a judicious application of Kasteleyn's representation and a tricky asymptotic analysis of Toeplitz matrices enabled Mc Coy and Wu \cite{\McCoyWu} to derive the precise asymptotic behavior of the two point function between the spins of two distant sites. 
From these computations, it is possible to give via \eqref{\eqdefsurften} a beautiful formula  for the surface tension.   
Their result states that for all $\beta>\beta_c$, and $x=(x_1, x_2)\in\Z^2$ : 
$$
\T(x)=x_1\ac(\sqrt{1+s^2x_1^2})+x_2\ac(\sqrt{1+s^2x^2_2}),\eqnum
$$\eqlabel{\formulasurften}   
where $s$ solves 
$$
\sqrt{1+s^2x_1^2}+\sqrt{1+s^2x_2^2}=\sinh(2\beta)+{1\over\sinh(2\beta)}.
$$
This formula is the only input from explicit computations that we use in this paper.
In \cite{\McCoyWu}, the formulas describing the asymptotic behavior of the correlation function and the surface tension are not given in  the form above. 
In the last section of this paper, we will show how the results of Mc Coy and Wu can be written in the form \eqref{\formulasurften}. As we will see, this rewriting reveals a simple connection between the geometry of the 2d-Ising Wulff shape and the large deviation rate function of the simple random walk on $\Z^2$. 

\noi The formula \eqref{\formulasurften} is the result of the limit \eqref{\eqdefsurften}, i.e., $n\rightarrow\infty$ but at a fixed $\beta$. 
The main purpose of our work is to derive the asymptotic behavior of the quantity in \eqref{\eqdefsurften} in a joint limit $\beta\da\beta_c$ and $n\rightarrow\infty$. Let us state the main result of this paper. 
\subhead{\subheadnum The main result}\endsubhead
\proclaim{Theorem \procnum}For all $x\in\Z^2$ and for any double sequence $n\ua\infty$ and $\beta\da\beta_c$ satisfying 
$$
{n\over\log n}(\beta-\beta_c)\ra\infty,\eqnum
$$\eqlabel{\regime}

\noi we have that
$$
-\lim_{n,\beta}{1\over (\beta-\beta_c)n}\log\mu^{\wh\beta}[\s(0)\s(nx)]=4|x|,\eqnum
$$\eqlabel{\mainresult}

\noi where $|x|$ is the ordinary Euclidean norm of $x$.
\endproclaim\proclabel{\jointlimit}
\noi In words, the above result states that in a joint limit satisfying \eqref{\regime}, it is still possible to define a rescaled surface tension and that the norm associated to this surface tension is the ordinary Euclidean norm.
The latter fact is an indication that the model studied in the regime \eqref{\regime} is asymptotically rotation-invariant.  
The regime of Theorem \procref{\jointlimit} is, up to a logarithmic correction, sharp. 
Indeed, in \cite{\BARROUCHMcCOYTRACYWU, \McCOYTRACY, \TRACY} it has been shown that in a regime where $n(\beta-\beta_c)=t$ stays constant, one gets that 
$$
\lim_{n,\beta}{\mu^{\beta}[\s(0)\s(ne_1)]\over n^{1/4}}=\eta(t),\eqnum
$$\eqlabel{\formulascalinglimit}

\noi where $e_1$ is the unit vector in the horizontal direction and $\eta(t)$ satisfies a Painlev\'e equation.
Thus, Theorem \procref{\jointlimit} delimits the threshold between the constant temperature case and the regime considered in \eqref{\formulascalinglimit}. The existence of the joint limit \eqref{\mainresult} is proved in two steps. First, we use \eqref{\formulasurften} to compute the limit  \eqref{\mainresult} when we first take the limit $n\uparrow\infty$ and then the limit $\beta\da\beta_c$. In the second step, which is our main contribution, we provide a probabilistic argument that shows that in any joint regime satisfying \eqref{\regime} the joint limit of \eqref{\mainresult} is well defined and equals the limit obtained in the first step. 
First, we give the easy step of the proof of Theorem \procref{\jointlimit}.
\proclaim{Proposition \procnum} Uniformly over $x\in\R^2$, we have that 
$$
\lim_{\beta\da\beta_c}\lim_{n\uparrow\infty}-{\log\mu^{\wh\beta}[\s(0)\s(\lfloor nx\rfloor)]\over (\beta-\beta_c)n|x|}=\lim_{\beta\da \beta_c}-{\T(x)\over (\beta-\beta_c)|x|}=4,
$$
where for $x=(x_1, x_2)$, $\lfloor nx\rfloor=(\lfloor nx_1\rfloor, \lfloor nx_2\rfloor)$ and for $i\in\{1, 2\}$ $\lfloor nx_i\rfloor$ is the largest integer which is smaller than or equal to $nx_i$.
\endproclaim\proclabel{\convEucNorm}

\noi The proof of Proposition \procref{\convEucNorm} is obtained directly from the formula \eqref{\formulasurften} by taking the limit $\beta\da \beta_c$.
\subhead\subheadnum Comparison with independent site percolation\endsubhead
Recently, a lot of progress has been made in the study of criticality in the context of planar independent site percolation on the triangular lattice $\Bbb T$. 
In this section, we compare the known results on the Ising model with their analogue in independent site percolation. The percolation analogue of  $\T$  is given by 
$$
\forall p>1/2\q\forall x\in\Bbb T\q \tau_p(x)=-\lim_{n\ra\infty}{1\over n}\log P^{\wh p}(0\lra nx),
$$ 
where $\wh p=1-p<1/2$ and where $P^{\wh p}$ is the probability measure corresponding to the site percolation of parameter $\wh p$ on $\Bbb T$.

\noi Proposition \procref{\convEucNorm} is actually a very strong statement that  implies the existence of the correlation length exponent $\nu=1$ in a rather strong from. It shows that the correlation length, if correctly rescaled becomes isotropic near the critical point. 
This result is at present time impossible to obtain without explicit computations. In percolation, Smirnov and Werner \cite{\SMIRNOVWERNER} used the convergence of critical percolation to the Schramm Loewner Evolution \cite{\SMIRNOV, \CAMIANEWMAN} and the scaling relations of Kesten \cite{\KESTEN} to prove that
$$
\forall x\in\Bbb T\q\lim_{p\da1/2}{\log\tau_p(x)\over\log(p-1/2)}=4/3.
$$
And the question whether the limit 
$
\lim_{p\da1/2}(p-1/2)^{-4/3}\tau_p(x)
$
exists is still open.

\noi The appearance in \eqref{\formulascalinglimit} of a function $\eta(t)$ which is related to a Painlev\'e equation is very striking. Such asymptotics are unavailable for percolation. 
However, an adequate modification of the critical scaling limit has been proposed in \cite{\CamiaNewmanII} to investigate a regime  analogous to the one considered in \eqref{\formulascalinglimit}.

\subhead{\subheadnum Organization of the paper}\endsubhead
In the next section we introduce the necessary notations and prove Theorem \procref{\jointlimit}. 
The third section contains a rewriting of the results \eqref{\formulasurften} of Mc Coy and Wu that permit to make the link with random walk and to draw some heuristics that explain the appearance of isotropy near the critical point. 

\head{\headnum Proof of the Theorem}\endhead
\noi The proof relies basically on sub-additivity. 
The easier part is the upper bound:
\proclaim{Lemma \procnum}For all $x\in\Z^2$ and for all double sequences $n\ua\infty$ and $\beta\da \beta_c$, we have that 
$$
\limsup_{n,\beta}{1\over{n(\beta-\beta_c)}}\log\mu^{\wh\beta}[\s(0)\s(nx)]\le-4|x|.
$$
\endproclaim\proclabel{\upperbound}
\demo{Proof}
Let us fix $\wh\beta<\beta_c$ and $x\in\Z^2$. By the FKG-inequality, the sequence 
$$
\left(\log\mu^{\wh\beta}[\s(0)\s(nx)],n\ge 1\right),
$$
is super-additive and thus, for every fixed $n>0$ and $\wh\beta<\beta_c$, 
$$
\mu^{\wh\beta}[\s(0)\s(nx)]\leq\exp\left(\lim_{m\ua\infty}{1\over m}\log\mu^{\wh\beta}[\s(0)\s(mnx)]\right)=\exp\left(-n\T(x)\right).\eqnum
$$\eqlabel{\expdecay}

\noi The result follows from Proposition \procref{\convEucNorm} that guarantees that 
$$
\lim_{\beta\da \beta_c}{1\over \beta-\beta_c}\T(x)=4|x|.
$$
\qed
\enddemo
\noi To prove the lower bound, we consider a sub-additive quantity that approximates adequately the two point function.  
To construct this quantity, we introduce for each $\beta>\beta_c$, the unit ball of the norm defined by $\T$: 
$$
U^{\beta}=\{y\in\R^2:\,\T(y)\leq 1\}.
$$
\par
\medskip%
\midinsert%
\vbox{\line{\hfil\epsfxsize=12cm \epsfbox{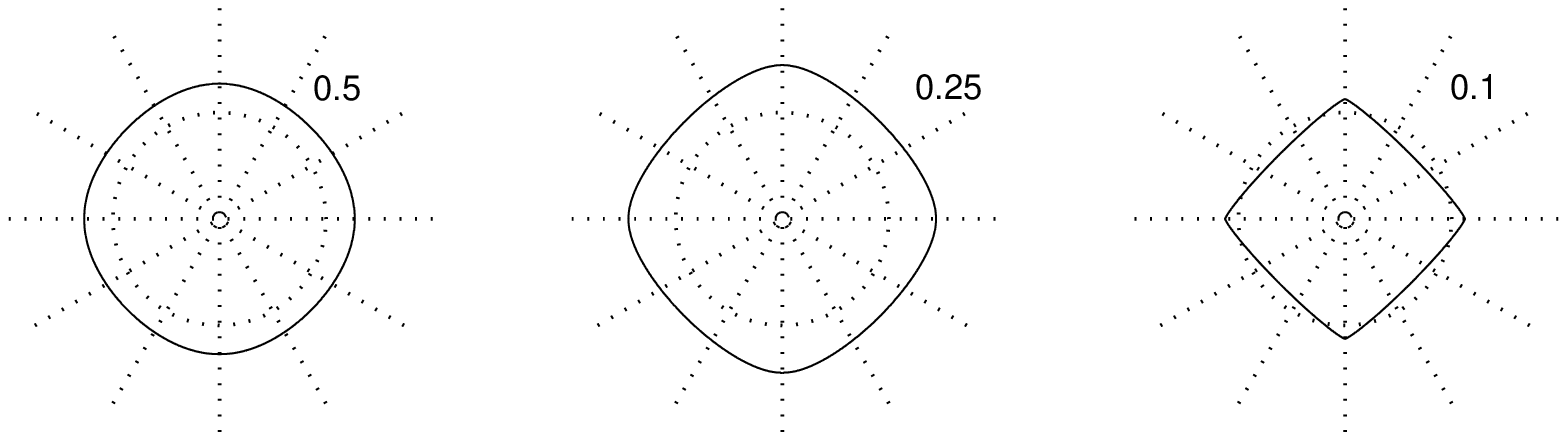}\hfil}\medskip}
\centerline{$U^{\beta}$ for increasing $\beta$, the number represents the radius of the dashed circle.}
\endinsert
\noi Next, for fixed $\beta>\beta_c$, $n>0$ and $x\in\Z^2$, we define the tangent line $w^{\beta}(nx)$ of $\T(nx)U^{\beta}$ at the point $nx$. We will also denote by 
$w^{\beta,+}(nx)$ the half plane delimited by $w^{\beta}(nx)$ and which does not contain the origin. 
We will also need the following discrete approximation of $w^{\beta}(nx)$:
$$
W^{\beta}(nx)=\{y\in w^{\beta,+}(nx)\cap\Z^2:\,0\leq(y-nx)\cdot v^{\beta}(x)\leq 1\},\eqnum
$$\eqlabel{\defW}%
where $v^{\beta}(x)$ is the unit length vector which is perpendicular to $w^{\beta}(nx)$ and which points towards the interior of $w^{\beta,+}(nx)$.

\noi\vbox{\centerline{
\psset{unit=0.08cm, linewidth=0.2}
\pspicture(-5,-10)(45,65)
\psline{->}(0,-5)(0,60)
\psline{->}(-5,0)(80,0)
\pspolygon[fillcolor=gray, fillstyle=vlines, hatchwidth=0.1, linecolor=white](-6, 50.5)(1.5, 60.5)(71.5, 8)(64, -2)
\pscurve(0,40)(6,39.5)(10,38)(26,26)
\pscurve(40,0)(39.5,6)(38,10)(26,26)
\psdot(19, 31.75)
\uput{2}[-90](19, 30.75){$nx$}
\psline(-6, 50.5)(64, -2)
\psline(1.5, 60.5)(71.5,8)
\psline[linewidth=0.4](-6, 50.5)(64, -2)
\psline{->}(19,31.75)(34, 51.75)
\uput{2}[90](34, 51.75){$v^{\beta}(x)$}
\uput{2}[180](-8, 50.5){$w^{\beta}(nx)$}
\psline[linewidth=0.1]{->}(60,35)(35, 25)
\uput{2}[0](60, 35){$W^\beta(nx)$}
\psline[linewidth=0.1]{->}(-10, 20)(38,10)
\uput{2}[180](-10, 20){$\partial U^\beta$}
\psline[linewidth=0.4](20.875, 34.25)(23.375, 32.375)(21.5,29.825)
\psline[linestyle=dashed](0,0)(19,31.75)
\endpspicture
}
}

\noi The quantity that permits us to establish the lower bound is 
$$
g(nx)=\sum_{y\in W^{\beta}(nx)}\mu^{\wh\beta}[\s(0)\s(y)].\eqnum
$$\eqlabel{\defg}

\noi The study of $g(nx)$ started in \cite{\CHAYESCHAYESCAMPANINO} and in \cite{\ALEXANDER}. 
The results of the last two cited papers where then refined \cite{\IOFFE, \CAMPANINOIOFFEVELENIKI, \CAMPANINOIOFFEVELENIKII, \CAMPANINOIOFFEVELENIKIII} to obtain 
Ornstein-Zernike asymptotics for the two point function.
Even though we establish rougher estimates than Ornstein-Zernike asymptotics, the quantity $g(nx)$ turns out to be important when analyzing the two point function near criticality.  
\proclaim{Lemma \procnum}For each fixed $\beta>\beta_c$ and $x\in\Z^2$, we have that:
\medskip
{\parindent 0.5cm
\item{i)}The following limit exists
$$
\lim_{n\ua\infty}-{1\over n}\log\sum_{y\in W^{\beta}(nx)}\mu^{\wh\beta}[\s(0)\s(y)],
$$ 
and is equal to the surface tension $\T(x)$.
\item{ii)}There exists a positive constant $K$ which does not depend on $x, n$ nor in $\beta$ such that
$$
\forall n\geq 1\qq\sum_{y\in W^{\beta}(nx)}\mu^{\wh\beta}[\s(0)\s(y)]\geq\exp(-n\T(x)-K).
$$
}
\endproclaim\proclabel{\hitline}
\demo{Proof}
The proof is an application of Simon's correlation inequality \cite{\SIMON} which states that for any two sites $x,y\in\Z^2$ and for any set $W$ separating the two sites 
$x,y$, the following holds
$$
\mu^{\wh\beta}[\s(x)\s(y)]\leq\sum_{z\in W}\mu^{\wh\beta}[\s(x)\s(z)]\mu^{\wh\beta}[\s(z)\s(y)].\eqnum
$$\eqlabel{\corrineqsimon}

\noi Note that for independent percolation a similar inequality holds thanks to the van den Berg-Kesten inequality. 

\noi Fix $n,m>0$ and consider
$$
g((n+m)x)=\sum_{y\in W^{\beta}((n+m)x)}\mu^{\wh\beta}[\s(0)\s(y)].
$$
For every $y\in W^{\beta}((n+m)x)$, the set $W^{\beta}(nx)$ separates the site $0$ from $y$, so that by \eqref{\corrineqsimon} we obtain
$$
g((n+m)x)\leq\sum_{y\in W^{\beta}((n+m)x)}\sum_{z\in W^{\beta}(nx)}\mu^{\wh\beta}[\s(0)\s(z)]\mu^{\wh\beta}[\s(z)\s(y)].
$$
By translation invariance of the measure $\mu^{\wh\beta}$  we get
$$
g((n+m)x)\leq\sum_{z\in W^{\beta}(nx)}\sum_{y\in W^{\beta}((n+m)x)}\mu^{\wh\beta}[\s(0)\s(z)]\mu^{\wh\beta}[\s(0)\s(y-z)].
$$
By \eqref{\defW}, for every $z\in W^\beta(nx)$ and $y\in W^\beta((n+m)x)$, we have that 
$$
0\leq(y-z-mx)\cdot v^\beta(x)\leq 2.
$$
Therefore, there exists a site $s\in W^\beta(mx)$ such that $|s-(y-z)|\leq 2$. 
Using the the FK-representation Theorem 1.16 of \cite{\GRIM} and the finite energy property Theorem 4.17 of \cite{\GRIM}, we get that 
$$
\mu^{\wh\beta}[\s(0)\s(y-z)]\leq e^{K}\mu^{\wh\beta}[\s(0)\s(s)],
$$  
where $K$ is a positive constant independent of $\beta$.
Hence
$$\ea{
g((n+m)x)\leq&e^K\sum_{z\in W^{\beta}(nx)}\sum_{s\in W^{\beta}(mx)}\mu^{\wh\beta}[\s(0)\s(z)]\mu^{\wh\beta}[\s(0)\s(s)]\cr
\leq&e^Kg(nx)g(mx).
}$$
By standard sub-additivity arguments, we get that the limit 
$$
\lim_{n\ua\infty}{1\over n}\log g(nx)=-g_\beta(x)
$$  
exists and moreover, for every $n\geq 1$
$$
g(nx)\geq\exp(-ng_\beta(x)-K).
$$
So it remains to prove that $g_{\beta}(x)=\T(x)$. 
We write
$$
\sum_{y\in W^{\beta}(nx)}
\mu^{\wh\beta}
[\s(0)\s(y)]
=
\mu^{\wh\beta}
[\s(0)\s(nx)]
+
\sum_{y\in W^{\beta}(nx)\setminus\{nx\}}
\mu^{\wh\beta}
[\s(0)\s(y)]
$$
Using \eqref{\expdecay} and the definition of $W^{\beta}(nx)$, we have 
$$\ea{
\sum_{y\in W^{\beta}(nx)\setminus\{nx\}}
\mu^{\wh\beta}
[\s(0)\s(y)]
\leq&
\sum_{y\in W^{\beta}(nx)\setminus\{nx\}}
\exp(-\T(y))\cr
\leq&
\sum_{k=0}^{\infty}
\sum_{y\in W^{\beta}(nx)\atop k<\T(y)-\T(nx)\leq k+1}\exp(-\T(y)).\cr
}
$$
To evaluate the last sum, we have to establish an upper bound for the cardinality of the set 
$$
\left\{y\in W^{\beta}(nx):\,k<\T(y)-\T(nx)\leq k+1\right\},
$$ 
where $k\geq 0$. For this, we note that there exists two positive constants $c_1$ and $c_2$ such that 
$$
\forall y\in\Z^2\q c_1|y|\leq\T(y)\leq c_2|y|.
$$
From there we conclude that 
$$\ea{
&\left\{y\in W^{\beta}(nx):\,k<\T(y)-\T(nx)\leq k+1\right\}\subset\cr
\subset&\left\{\,y\in W^{\beta}(nx):\,|y|\in(\,c_2^{-1}(k+\T(nx)), c_1^{-1}(k+1+\T(nx))\,]\,\right\}.
}$$
Using the fact that $\T(nx)\leq c_2 n|x|$ we can bound the width of the interval containing $|y|$:
$$
c_1^{-1}(k+1+\T(nx))-c_2^{-1}(k+\T(nx))\leq {c_2-c_1\over c_1c_2}(k+c_2n|x|)+{1\over c_1}.
$$  
Therefore, for $n$ large enough we can find a positive constant $C$ such that 
$$
\left|\left\{y\in W^{\beta}(nx):\,k<\T(y)-\T(nx)\leq k+1\right\}\right|\leq C(k+1+n).
$$
Hence
$$\ea{
\sum_{y\in W^{\beta}(nx)\setminus\{nx\}}
\mu^{\wh\beta}
[\s(0)\s(y)]
\leq&C\exp(-n\T(x))\sum_{k=0}^\infty(k+1+n)e^{-k}\cr
\leq&C\left({ne\over e-1}+{e^2\over(e-1)^2}\right)\exp(-n\T(x)).
}
$$ 
Thus
$$
\lim_{n\ua\infty}
{1\over n}
\log
\sum_{y\in W^{\beta}(nx)}
\mu^{\wh\beta}
[\s(0)\s(y)]
=
\lim_{n\ua\infty}
{1\over n}
\log
\mu^{\wh\beta}
[\s(0)\s(nx)]
=
-\T(x).
$$
\qed
\enddemo
\noi It is crucial to note that in the case $x=ae_i, i\in\{1, 2\}$ where $a\in\N\setminus\{0\}$ and $(e_1, e_2)$ is a orthonormal basis of $\R^2$, then one can derive a slightly stronger result. 
\proclaim{Lemma \procnum}For each fixed $a\in\N\setminus\{0\}$, $\beta>\beta_c$ and $x=ae_i, i=1, 2$, we have that:
\medskip
{\parindent 0.5cm
\item{i)}The following limit exists
$$
\lim_{n\ua\infty}-{1\over n}\log\sum_{y\in w^{\beta}(nx)}\mu^{\wh\beta}[\s(0)\s(y)],
$$ 
and is equal to the surface tension $\T(x)$.
\item{ii)}For all $n>0$
$$
\sum_{y\in w^{\beta}(nx)}\mu^{\wh\beta}[\s(0)\s(y)]\geq\exp(-n\T(x)).
$$
}
\endproclaim\proclabel{\hitlineaxis}
\demo{Proof}
The proof is almost the same than the proof of Lemma \procref{\hitline}, except that $W^\beta(nx)$ can be replaced by the line $w^\beta(nx)$ then by translation invariance we have that for every integers $n$ and $m$ and for each $z\in w^\beta(nx)$ and $y\in w^\beta((m+n)x)$,  we have that 
$
y-z\in w^{\beta}(mx).
$ 
This turns $(g(nx), n\geq1)$ into an exact sub-additive sequence. This property permits to prove the enhancement ii) for such a particular $x$. 
\qed\enddemo
\noi The last two lemma indicate that \eqref{\defg} is a nice quantity to get the desired lower bound part of Theorem \procref{\jointlimit}. 
 It turns out that if $x$ points in the direction of one of the two vectors $(e_1, e_2)$ of the canonical orthonormal basis of $\Z^2$, then a judicious symmetry argument of Chayes, Chayes and Campanino \cite{\CHAYESCHAYESCAMPANINO} permits to obtain the desired lower bound.
\proclaim{Proposition \procnum}
There exists $\beta_0>\beta_c$ and  a positive constant $c$, such that for all $\beta_c<\beta<\beta_0$, for all $n>1$ and for all $x=ae_i, i=1,2$ where $a\in\N\setminus\{0\}$ the following holds

$$
\mu^{\wh\beta}[\s(0)\s(nx)]\geq\left({1-cne^{-(\beta-\beta_c)n/2}{2-e^{-(\beta-\beta_c)n/2}\over(1-e^{-(\beta-\beta_c)n/2})^2}\over(6|x|+1)n}\right)^2 \exp(-n\T(x)),
$$
so that for every double sequence $n\ua\infty$ and $\beta\da \beta_c$ satisfying {\rm\eqref{\regime}}, we have 
$$
\lim_{n, \beta}
{1\over (\beta-\beta_c)n}\log\mu^{\wh\beta}[\s(0)\s(nx)]=-4|x|.
$$
\endproclaim\proclabel{\lowerboundaxis}
\demo{Proof}
Let $n>0$ and $\beta>\beta_c$. By Lemma \procref{\hitlineaxis}, we have that
$$ 
\sum_{y\in w^{\beta}(nx)}\mu^{\wh\beta}[\s(0)\s(y)]\geq\exp(-n\T(x)).
$$
On the other side, 
$$\ea{
&\sum_{y\in w^{\beta}(nx)\atop \T(y)-n\T(x)>(\beta-\beta_c)n}
\mu^{\wh\beta}[\s(0)\s(y)]\leq\cr
\leq&e^{-n\T(x)}\sum_{k=1}^{+\infty}|\{y\in w^{\beta}(nx):\,k<{\T(y)-\T(nx)\over{(\beta-\beta_c)n}}\leq k+1\}|\,e^{-k(\beta-\beta_c)n}.
}$$
By Proposition \procref{\convEucNorm}, there exists $\beta_0>\beta_c$ such that 
$$
\forall \beta_c<\beta<\beta_0\,\forall z\in\Z^2\q2|z|\leq{\T(z)\over\beta-\beta_c}\leq6|z|.\eqnum
$$\eqlabel{\tauLtwo}%
Therefore for all $k>0$
$$\ea{
&\left\{y\in w^{\beta}(nx):\,k<{\T(y)-\T(nx)\over{(\beta-\beta_c)n}}\leq k+1\right\}\subset\cr
\subset&\left\{y\in w^{\beta}(nx):\,|y|\in(\,{kn\over 6}+{n\over 3}|x|,\, {(k+1)n\over 2}+3n|x|\,]\right\}.
}$$
From where we conclude that there exists a positive constant $c=c(x)$, such that 
$$
|\{y\in w^{\beta}(nx):\,k<{\T(y)-\T(nx)\over{(\beta-\beta_c)n}}\leq k+1\}|\leq cn(k+1).
$$
And hence
$$\ea{
\sum_{y\in w^{\beta}(nx)\atop \T(y)-n\T(x)>(\beta-\beta_c)n}\mu^{\wh\beta}[\s(0)\s(y)]
&\leq cne^{-n\T(x)}\sum_{k=1}^{+\infty}(k+1)\,e^{-k(\beta-\beta_c)n}\cr
\leq&cne^{-(\beta-\beta_c)n}{2-e^{-(\beta-\beta_c)n}\over (1-e^{-(\beta-\beta_c)n})^2}\exp(-n\T(x)).
}$$

\noi Thus 
$$
\sum_{y\in w^{\beta}(nx)\atop \T(y)-n\T(x)<(\beta-\beta_c)n}\hskip-0.5cm\mu^{\wh\beta}[\s(0)\s(y)]\geq\left[1-cne^{-(\beta-\beta_c)n}{2-e^{-(\beta-\beta_c)n}\over (1-e^{-(\beta-\beta_c)n})^2}\right] \exp(-n\T(x)).
$$
Next, we bound the cardinality of the summation set above. Using \eqref{\tauLtwo}, we have
$$\ea{
&\left\{y\in w^\beta(nx):\,\T(y)-n\T(x)<(\beta-\beta_c)n\right\}\subset\cr
\subset&\left\{y\in w^\beta(nx):\,|y|\leq(3|x|+{1/2)n}\right\}.
}$$
Therefore $|\{y\in w^\beta(nx):\,\T(y)-n\T(x)<(\beta-\beta_c)n\}|\leq(6|x|+1)n$. Thus there exists a site $y_{n,\beta}$ such that $\T(y_{n,\beta})-n\T(x)\leq (\beta-\beta_c)n$ and 
$$
\mu^{\wh\beta}[\s(0)\s(y_{n, \beta})]\geq{1-cne^{-(\beta-\beta_c)n}{2-e^{-(\beta-\beta_c)n}\over (1-e^{-(\beta-\beta_c)n})^2}\over(6|x|+1)n}\exp(-n\T(x)).\eqnum
$$ \eqlabel{\lbI}

\noi By symmetry with respect to $w^\beta(nx)$, for every $y_{n, \beta}\in w^{\beta}(nx)$:
$$
\mu^{\wh\beta}[\s(y_{n, \beta})\s(2nx)]=\mu^{\wh\beta}[\s(0)\s(y_{n, \beta})].
$$
Combining the last equality with the FKG inequality, we get that 
$$\ea{
\mu^{\wh\beta}
[\s(0)\s(2nx)]
\geq&
\mu^{\wh\beta}
[\s(0)\s(y_{n, \beta})]
\mu^{\wh\beta}
[\s(y_{n, \beta})\s(2nx)]\cr
\geq&
\mu^{\wh\beta}
[\s(0)\s(y_{n, \beta})]^2.
}$$
Thus, we get from \eqref{\lbI} that 
$$\ea{
\mu^{\wh\beta}
[\s(0)\s(2nx)]
\geq&\mu^{\wh\beta}[\s(0)\s(y_{n, \beta})]^2\cr
\geq&\left({1-cne^{-(\beta-\beta_c)n}{2-e^{-(\beta-\beta_c)n}\over (1-e^{-(\beta-\beta_c)n})^2}\over(6|x|+1)n}\right)^2\exp(-2n\T(x)).
}$$
Combining the last inequality with \eqref{\expdecay}, we deduce that for any joint limit $n\ua\infty$ and $\beta\da \beta_c$ that satisfies \eqref{\regime} we have that 
$$
\lim_{n,\beta}
{1\over (\beta-\beta_c)n}
\log
\mu^{\wh\beta}[\s(0)\s(nx)]
=
-4|x|.
$$
\qed
\enddemo
\noi Now we are ready to proceed to the completion of the proof of Theorem \procref{\jointlimit}.
\demo{Proof of Theorem \procref{\jointlimit}}

\noi Let $x\in\Z^2$ be fixed. The upper bound is contained in Lemma \procref{\upperbound}.
Thus, it remains to determine the conditions on the regime that guarantee that  
$$
\liminf_{n,\beta}
{1\over (\beta-\beta_c)n}
\log
\mu^{\wh\beta}
[\s(0)\s(nx)]
\geq
-4|x|
.
$$
From ii) of Lemma \procref{\hitline}, There exists a positive constant $K$ such that for all $n>0$ and $\beta>\beta_c$ 
$$
\sum_{y\in W^{\beta}(nx)}
\mu^{\wh\beta}
[\s(0)\s(y)]
\geq
\exp(-n\T(x)-K)
.
$$
Fix $\eps>0$. By Proposition \procref{\convEucNorm}, there exists $\beta(\eps)>\beta_c$ such that for all $\beta_c<\beta<\beta(\eps)$
$$
\forall\,z\in\Z^2\qq{\T(z)\over|z|}\in(4-\eps, 4+\eps).
$$ 
Therefore, for all $0<\eps<4$ and for all $\beta_c<\beta<\beta(\eps)$:
$$\ea{
|y|-|nx|>\eps n\q\Rightarrow\q&{1\over 4-\eps}(\T(y)-n\T(x))>{2\eps\over16-\eps^2}n\T(x)+\eps n(\beta-\beta_c)n\cr
\q\Rightarrow\q&\T(y)-n\T(x)>\eps(\beta-\beta_c)n.
}
$$
Using \eqref{\expdecay}, we get  
$$\ea{
&\sum_{y\in W^{\beta}(nx)\atop|y|-|nx|>\eps n}\mu^{\wh\beta}[\s(0)\s(y)]\leq\sum_{y\in W^{\beta}(nx)\atop \T(y)-\T(nx)>\eps(\beta-\beta_c)n}\mu^{\wh\beta}[\s(0)\s(y)]\cr
\leq&e^{-n\T(x)}\sum_{k=1}^{+\infty}|\{y\in W^{\beta}(nx):\, k<{\T(y)-\T(nx)\over\eps(\beta-\beta_c)n}\leq k+1\}|e^{-\eps(\beta-\beta_c)nk}.\cr
}$$
By the same arguments than those used in the proof of Proposition \procref{\convEucNorm}, we have that
$$
\left|\left\{y\in W^{\beta}(nx):\,k<{\T(y)-\T(nx)\over\eps(\beta-\beta_c)n}\leq k+1\right\}\right|\leq c_1n(k+1),
$$
where $c_1=c_1(x)$ is a positive constant. 
Hence 
$$\ea{
&\sum_{y\in W^{\beta}(nx)\atop|y|-|nx|>\eps n}\mu^{\wh\beta}[\s(0)\s(y)]\leq c_1ne^{-n\T(x)}\sum_{k=1}^\infty(k+1)e^{-\eps(\beta-\beta_c)nk}\cr
\leq&c_1ne^{-\eps(\beta-\beta_c)n}{2-e^{-\eps(\beta-\beta_c)n}\over(1-e^{-\eps(\beta-\beta_c)n})^2}\exp(-n\T(x)).
}$$
Therefore 
$$
\sum_{y\in W^{\beta}(nx)\atop|y|-|nx|\leq\eps n}\hskip-0.4cm\mu^{\wh\beta}[\s(0)\s(y)]\geq\left(e^{-K}-c_1ne^{-\eps(\beta-\beta_c)n}{2-e^{-\eps(\beta-\beta_c)n}\over(1-e^{-\eps(\beta-\beta_c)n})^2}\right)\exp(-n\T(x)).
$$
And for any joint limit satisfying \eqref{\regime}, we have  
$$
\liminf_{n,\beta}
{1\over (\beta-\beta_c)n}
\log
\sum_{y\in W^{\beta}(nx)\atop |y|-|nx|\leq\eps n}
\mu^{\wh\beta}
[\s(0)\s(y)]
\geq
-4|x|.
$$
Since $W^\beta(nx)$ is a cylinder of finite basis, there exists a positive constant $c_2=c_2(x)$ such that 
$$
\left|\{y\in W^\beta(nx):\,|y|-|nx|\leq\eps n\}\right|\leq c_2n.
$$
Therefore there exists a site $y_{\eps,n, \beta}$ in $W^{\beta}(nx)$ such that $|y|-|nx|\leq\eps n$ and 
$$
\mu^{\wh\beta}
[\s(0)\s(y_{\eps,n, \beta})]
\geq
{1\over c_2 n}\left(e^{-K}-c_1ne^{-\eps(\beta-\beta_c)n}{2-e^{-\eps(\beta-\beta_c)n}\over(1-e^{-\eps(\beta-\beta_c)n})^2}\right)\exp(-n\T(x)).
$$
So that for any joint limit satisfying  \eqref{\regime}
$$
\liminf_{n,\beta}
{1\over (\beta-\beta_c)n}
\log
\mu^{\wh\beta}
[\s(0)\s(y_{\eps,n, \beta})]
\geq
-4|x|.\eqnum
$$\eqlabel{\borninfI}

\noi On the other hand, by \eqref{\expdecay} and by the definition of $W^{\beta}(nx)$  we have that 
$$
\mu^{\wh\beta}[\s(0)\s(y_{\eps, n, \beta})]\leq\exp(-\T(y_{\eps, n, \beta}))\leq\exp(-n\T(x)).
$$
Combining the last inequality with \eqref{\borninfI} and using Proposition \procref{\convEucNorm}, we get that for any joint limit satisfying \eqref{\regime}, the following holds
$$
-\lim_{n,\beta}
{\log\mu^{\wh\beta}[\s(0)\s(y_{\eps, n, \beta})]\over (\beta-\beta_c)n}
=
\lim_{n,\beta}
{\T(y_{\eps, n, \beta}/n)\over (\beta-\beta_c)}
=\lim_{\beta\da \beta_c}
{\T(x)\over (\beta-\beta_c)}=4|x|.\eqnum
$$\eqlabel{\convy} 

\noi Next, by the FKG inequality and by translation invariance, we have that 
$$\ea{
&\M[\s(0)\s(nx)]\geq\M[\s(0)\s(y_{\eps, n, \beta})]\M[\s(0)\s(nx-y_{\eps, n, \beta})]\cr
&\geq\M[\s(0)\s(y_{\eps, n, \beta})]\times\cr
&\times\M[\s(0)\s(((n(x-y_{\eps, n,\beta}/n))\cdot e_1)e_1)]\,\M[\s(0)\s(((n(x-y_{\eps, n, \beta}/n))\cdot e_2)e_2)].
}\eqnum$$\eqlabel{\borninfII}

\noi To finish the proof, we need the following lemma: 
\proclaim{Lemma \procnum}
Let $\eps>0$. For all $x\in\Z^2$ there exists $\beta(\eps, x)>\beta_c$ such that for all $\beta_c<\beta<\beta(\eps, x)$ and for all $n>1$ the following holds 
$$
\{y\in W^{\beta}(nx):\,|y/n|-|x|\leq\eps\}\subset\{y\in W^{\beta}(nx):\,|y/n-x|\leq\sqrt{(4x+2\eps)\eps}+{1\over n}\}.
$$
\endproclaim\proclabel{\normcontrol}

\noi We postpone the proof of Lemma \procref{\normcontrol} to the end of this section and continue the proof of the theorem. By Propoposition \procref{\lowerboundaxis} we have that
$$\ea{
\M[\s(0)\s(((n(x-y_{\eps, n,\beta}/n))\cdot e_i)e_i)]\geq&\left({1-cne^{-(\beta-\beta_c)n/2}{(2-e^{-(\beta-\beta_c)n/2})\over(1-e^{-(\beta-\beta_c)n/2})^2}\over(6|x|+1)n}\right)^2\times\cr
&\times\exp(-n\T(((x-y_{\eps, n, \beta}/n)\cdot e_i)e_i)).
}$$
By Lemma \procref{\normcontrol}, there exists $\beta(\eps, x)>\beta_c$ such that:
$$\forall\,\beta_c<\beta<\beta(\eps, x),\,n>1:\q((x-y_{\eps, n, \beta}/n)\cdot e_i)e_i|\leq\sqrt{(4x+2\eps)\eps}+{1\over n}.$$ 
Hence, for any joint limit satisfying \eqref{\regime}, we obtain that 
$$ 
\liminf_{n,\beta}
{1\over (\beta-\beta_c)n}
\log
\M[\s(0)\s(((n(x-y_{\eps, n, \beta}/n))\cdot e_i)e_i)]
\geq -4\sqrt{(4x+2\eps)\eps}.$$
Combining the last result with \eqref{\borninfII} and \eqref{\convy}, we obtain that for all $\eps>0$
$$
\liminf_{n,\beta}
{1\over (\beta-\beta_c)n}
\log
\M[\s(0)\s(nx)]\geq -4|x|-8\sqrt{(4x+2\eps)\eps}.
$$
Since the inequality is true for any $\eps>0$, we get the desired result. 
\qed
\enddemo
\demo{Proof of Lemma \procref{\normcontrol}}
The problem reduces to the following situation: 
given a circle $C(O,|x|+\eps)$ centered at the origin 0 and of radius $|x|+\eps$, find a circle $C(x,\delta(\eps, x))$ centered at $x$ of radius $\delta(\eps,x)$, such that $C(O,|x|+\eps)\cap w^{\beta}(x)\subset C(x,\delta(\eps, x))$. 
Let us denote by $\alpha(\beta, x)$ the angle formed by $w^{\beta}(x)$ and the line passing through $x$ and perpendicular to $(O,x)$, see the figure below. 

\noi\vbox{\centerline{
\psset{unit=0.08cm, linewidth=0.2}
\pspicture(0,-35)(20,70)
\rput{35}
{
\psdots(0,0)(25,0)
\psline(0,0)(25,0)
\uput[0]{-35}(24,0){$x$}
\uput[0]{-35}(-9,0){$O$}
\uput[0]{-35}(32,4){$\delta(\eps,x)$}
\uput[0]{-35}(-12,-13){$\eps+|x|$}
\uput[0]{-35}(8,56){$\alpha(\beta, x)$}
\uput[0]{-35}(35,-40){$w^{\beta}(x)$}
\psline(25,-50)(25,60)
\psline(25,0)(13,48)
\psline(25,0)(35,-40)
\psline(22,0)(22,-3)
\psline(22,-3)(25,-3)
\psline[linestyle=dashed]{<->}(25,0)(44,13)
\psline{<->}(0,0)(-21,-21)
\psarc(25,0){40}{90}{104}
\pscircle(0,0){30}
\pscircle[linestyle=dashed](25,0){23.5}
}
\endpspicture
}
}
Elementary geometric considerations show that $\delta$ can be chosen as follows 
$$
\delta(x,\eps)=\sqrt{2|x|^2\sin^2(\alpha)+2|x|\left(\eps+|\sin(\alpha)|\sqrt{\sin^2(\alpha)|x|^2+2\eps|x|+\eps^2}\right)+\eps^2}.
$$
By Proposition \procref{\convEucNorm}, we have that $\alpha(\beta, x)\rightarrow 0$, when $\beta\da \beta_c$. 
Thus, for every $\eps>0$ and for every $x\in\Z^2$, there exists $\beta(\eps, x)>\beta_c$ such that 
$$\forall\,\beta_c<\beta<\beta(\eps, x):\q\delta(x,\eps)<\sqrt{2(2\eps|x|+\eps^2)}.$$ 
Replacing $w^{\beta}$ by $W^{\beta}$, induces an extra error of $1/n$ and we are done.
\qed
\enddemo

\head{\headnum Asymptotics of correlation function and random walks}\endhead  
As promised in the introduction, we include in this section some explanations of the results of Mc Coy and Wu about asymptotics of correlations \cite{\McCoyWu}. 
On page 305 of \cite{\McCoyWu}, formula (4.38) and (4.39) give an asymptotic expansion of $\mu^{\wh\beta}[\s(0)\s(x)]$, with $x=(M,N)$, when $\wh\beta<\beta_c$ is hold fixed and 
$M^2+N^2\ua\infty$.
Actually, a closer look to the computations shows that (4.39) is just the expansion of the following double integral given by (4.22) and (4.23) of \cite{\McCoyWu}. 
So that the results of \cite{\McCoyWu} can be rewritten as
$$
\mu^{\wh\beta}
[\s(0)\s(x)]
\sim
{
[
\sinh^{-4}(2\wh\beta)
-1
]^{1/4}
\over
4\pi^2\gamma_\beta a_\beta
}
\int_{-\pi}^{\pi}
\int_{-\pi}^{\pi}
d\theta_1\,d\theta_2
{\cos(M\theta_1+N\theta_2)
\over 
1-{\gamma_\beta\over a_\beta}(\cos\theta_1+\cos\theta_2)
},\eqnum
$$\eqlabel{\formuleMcCoyWuI}

\noi when $\wh\beta<\beta_c$ is fixed and $|x|^2=M^2+N^2\ua\infty$. Where $a_\beta=(1+\tanh^2(\wh\beta))^2$ and $\gamma_\beta=2\tanh(\wh\beta)(1-\tanh^2(\wh\beta))$.
Instead of expanding this integral, we remark that this is just the generating function of a simple random walk.  
Indeed, let $(S_n)_{n\geq 1}$ be a simple symmetric random walk on $\Z^2$ starting at the origin and let $m\in(0,1)$, then it is known that for every $x=(x_1, x_2)\in\Z^2$:
$$
\sum_{k=0}^{\infty}P(S_k=x)m^k={1\over 4\pi^2}\int_{-\pi}^{\pi}\int_{-\pi}^{\pi}d\theta_1d\theta_2{\cos(x_1\theta_1+x_2\theta_2)\over{1-{m\over 2}(\cos(\theta_1)+\cos(\theta_2))}}.
$$
Hence, if we denote by $\xi_m$ a geometric random variable of parameter $m$ and independent of the walk then \eqref{\formuleMcCoyWuI} can be rewritten as 
follows 
$$
\mu^{\wh\beta}
[\s(0)\s(x)]
\sim
{
[
\sinh^{-4}(2\wh\beta)
-1
]^{1/4}
\over
\gamma_\beta a_\beta
}
E(V_{m(\beta)}(x)),\eqnum
$$\eqlabel{\formuleMcCoyWuII}

\noi here $V_m(x)=\sum_{k=0}^{\xi_m}1_{\{S_k=x\}}$ is the number of visits the killed random walk do to the site $x$. 
The surviving rate $m(\beta)$ of the walk is given by the following formula 
$$
m(\beta)={2\gamma_\beta\over a_\beta}={2\over{\sinh(2\wh\beta)+{1\over\sinh(2\wh\beta)}}}={2\over\sinh(2\wh\beta)+\sinh(2\beta)}.\eqnum
$$\eqlabel{\formulemT}

\par\noi
\vbox{\line{\hfil\epsfxsize=8cm \epsfbox{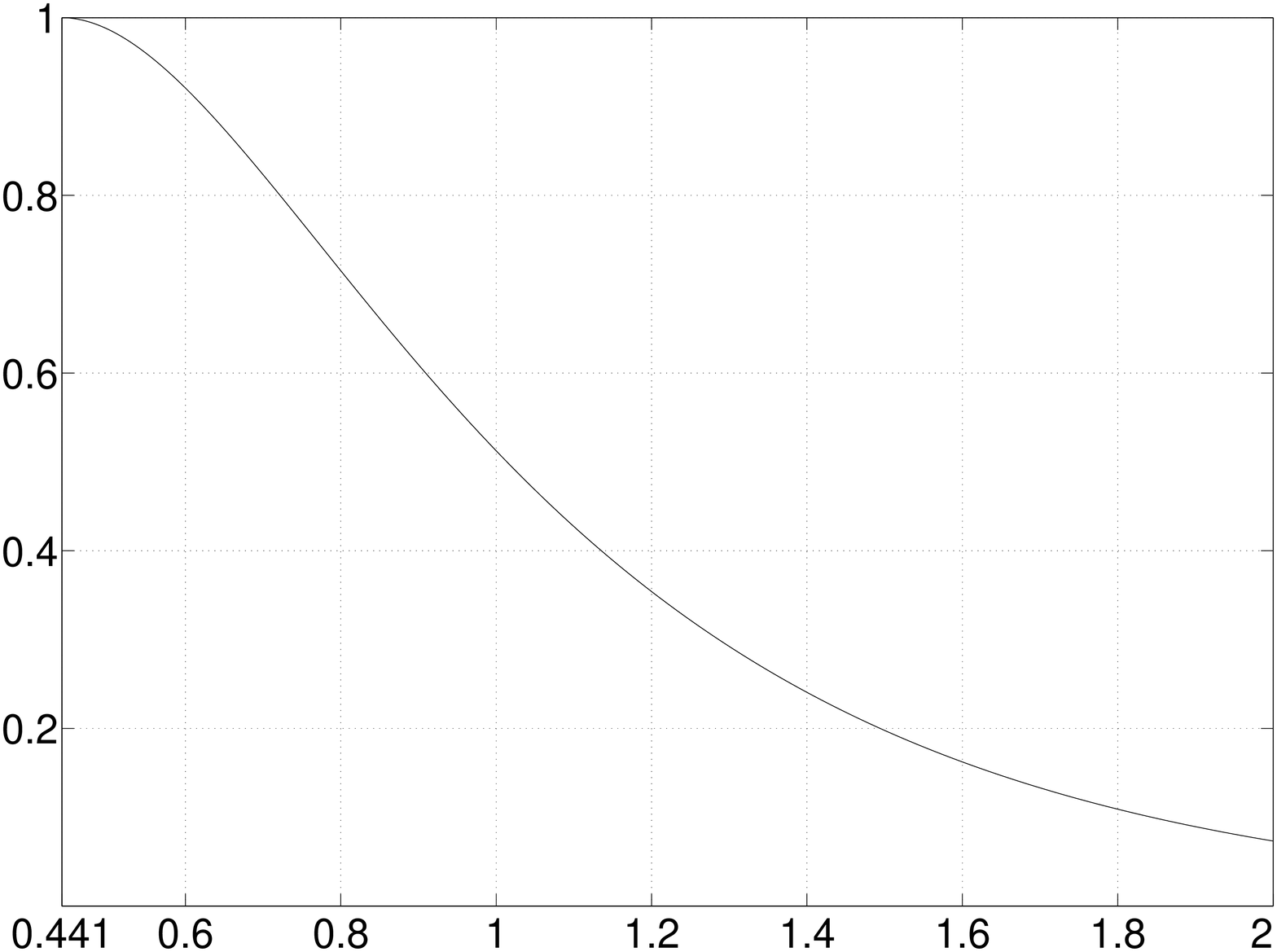}\hfil}\medskip}
\centerline{$m(\beta)$ as a function of $\beta_c\leq\beta\leq 2$.}
\bigskip
\noi In words, it is possible to calibrate the surviving rate of a killed random walk, such that the asymptotics of the two point function of the 2d-Ising model can be described in terms of the random walk. 

\subhead{\subheadnum The surface tension at fixed temperature}\endsubhead
In this section we show how to derive the exact formula for the surface tension from \eqref{\formuleMcCoyWuII} using large deviation estimates for random walks.
This will give a nice expression for $\T$ involving the Cram\'er function of the random walk. 
\proclaim{Proposition \procnum} For all $x\in\R^2$ and for all $\beta>\beta_c$, we have that
$$
\T(x)=\inf_{\gamma>0}\left(\lambda(\beta)\gamma+\gamma I(x/\gamma)\right),
$$ 
where $I(\cdot)$ is the Cram\'er function of the simple random walk on $\Z^2$ and 
$$
\lambda(\beta)=-\log(m(\beta))=\log\left({\sinh(2\wh\beta)+\sinh(2\beta)\over 2}\right).
$$ 
\endproclaim\proclabel{\tauCramer}
\demo{Proof}
It is sufficient to prove the result for $x\in\Z^2$. The general case follows from the fact that $\T$ is a norm. 
From \eqref{\eqdefsurften} and from \eqref{\formuleMcCoyWuII}, we know that
$$
\T(x)=-\lim_{n\ra\infty}{1\over n}\log E[V_{m(\beta)}(nx)],
$$
where $m(\beta)$ is given by \eqref{\formulemT}. 
For every $\beta>\beta_c$, we introduce the quantity $\lambda=\lambda(\beta)=-\log(m(\beta))$. Also, for any $n\geq 1$ and $x\in\Z^2$, we define the first time the random walk hits $nx$ 
$$
H(nx)=\inf\{k\geq0: S_k=nx\}.
$$
It is elementary to check that 
$$
E[\exp(-\lambda H(nx))]
\leq 
E[V_{m(\beta)}(nx)]
\leq
{1\over 1-m(\beta)}E[\exp(-\lambda H(nx))],
$$
So that 
$$
\T(x)=-\lim_{n\ra\infty}{1\over n}\log E[\exp(-\lambda H(nx))].
$$
We start with a lower bound. Let us fix $\gamma>0$, then  
$$\eqalign{
E[\exp(-\lambda H(nx))]
&
\geq E[\exp(-\lambda H(nx)),H(nx)\leq\gamma
n]\cr
&\geq \exp(-\lambda\gamma n)P(\exists k<\gamma n: S_{k}=nx).
}$$
Fix $\eps>0$. We notice that if there exists $k\geq 0$ such that  $S_k=nx$, then for all $y\in B(nx, \eps n)$ we have 
$$
P(S_k=nx)\geq\left({1\over 4}\right)^{\eps n}P(S_k=y). 
$$
Therefore, there exists a positive constant $c$ such that 
$$\ea{
\log P(\exists k<\gamma n: S_{k}=nx)\geq&-\eps n\log 4+\sup_{y\in B(nx, \eps n)}P(\exists k<\gamma n: S_k=ny)\cr
\geq&-\eps n\log 4-c\log n+P(\exists k<\gamma n: S_k\in B(nx, \eps n)),
}$$
where $B(nc, \eps n)=\{y\in\Z^2:\,|y-nx|\leq\eps n\}$.

\noi By further restricting $k$ to be equal to the largest integer $\lfloor \gamma n\rfloor$ which is smaller than $\gamma n$, we get that 
$$
\liminf_{n\ra\infty}{1\over n}\log E[\exp(-\lambda H(nx))]\geq-\eps\log 4+\limsup_{n\ra\infty}{1\over n}\log P(S_{\lfloor \gamma n\rfloor}\in B(nx, \eps n)).
$$
From Cram\'er's Theorem \cite{\DZ} for the simple random walk, we have 
$$
\limsup_{n\ra\infty}{1\over n}\log P(S_{\lfloor \gamma n\rfloor}\in B(nx, \eps n))\geq-\inf_{y\in B(x, \eps)}\gamma I(x/\gamma). 
$$
By the continuity of the rate function 
$$
I(x_1, x_2)=\sup_{(\lambda_1, \lambda_2)\in\R^2}\left(x_1\lambda_1+x_2\lambda_2-\log\left({1\over 2}(\cosh(\lambda_1)+\cosh(\lambda_2))\right)\right)
$$
we get that
$$
\liminf_{n\ra\infty}{1\over n}\log
E[\exp(-\lambda H(nx))]\geq -2\eps\log 4-\gamma I(x/\gamma)-\gamma\lambda.
$$
Since the result holds for all $\eps>0$ and for any $\gamma>0$, we get that 
$$
\liminf_{n\ra\infty}{1\over n}\log E[\exp(-\lambda H(nx))]\geq-\inf_{\gamma>0}\left(\lambda\gamma+\gamma I(x/\gamma)\right),
$$
It remains to establish the corresponding upper bound. For this, let us fix $\nu>0$ and consider
$$
E[\exp(-\lambda H(nx)),\, H(nx)\leq \nu n].
$$
For any large integer $M$, we write 
$$\ea{
&\log E[e^{-\lambda H(nx)},\,H(nx)\leq\nu n]\cr
\leq&\log\sum_{k=0}^{M-1}E[e^{-\lambda H(nx)},\frac{k}{M}\nu n<H(nx)\leq\frac{k+1}{M}\nu n)\cr
\leq&\log\sum_{k=0}^{M-1}\exp(-\lambda\frac{k}{M}\nu n)P(\frac{k}{M}\nu n<H(nx)\leq\frac{k+1}{M}\nu n).
}$$
Let $x=(x_1, x_2)$ and observe that $t<n(|x_1|+|x_2|)$ implies that $P(S_t=nx)=0$. 
Thus, we can restrict the above sum to those $k$ that are larger than or equal to $k_0=(M/\nu)(|x_1|+|x_2|)-1$. For $M$ large enough, $k_0$ is strictly positive. 
Fix such an $M$ and a $k\geq k_0>0$, we have 
$$
P({k\over M}\nu n<H(nx)\leq {k+1\over M}\nu n)\leq P(\exists t\in({k\over M}\nu n,\, {k+1\over M}\nu n]:\, S_t=nx).
$$
Notice that if there exists a $t\in ({k\over M}\nu n,\, {k+1\over M}\nu n]$ such that $S_t=nx$, then the random walk $\left(S_s, s\in({k\over M}\nu n,\, {k+1\over M}\nu n]\right)$ is included in $B(nx, \nu n/M)$. 
Therefore
$$\ea{
P({k\over M}\nu n<H(nx)\leq {k+1\over M}\nu n)\leq&P(\forall t\in ({k\over M}\nu n,\, {k+1\over M}\nu n]:\, S_t\in B(nx, {\nu n\over M}))\cr
\leq&P(S_{\lceil{k\over M}\nu n\rceil}\in nB(x, {\nu\over M})),
}$$
where $\lceil{k\over M}\nu n\rceil$ denotes the smallest integer greater than or equal to $\nu nk/M$.
By Cram\'er's Theorem, we get 
$$
\liminf_{n\ra\infty}{1\over n}\log P({k\over M}\nu n\leq H(nx)\leq {k+1\over M}\nu n)\leq -\inf_{y\in N(x,\nu/M)}{k\over M}\nu I\left({y\over {k\over M}\nu}\right).
$$
Using the continuity of the rate function $I$, we further get that for any $\eps>0$ there exists $M_0$ such that for any $M>M_0$
$$
\liminf_{n\ra\infty}{1\over n}\log P({k\over M}\nu n\leq H(nx)\leq {k+1\over M}\nu n)\leq -{k\over M}\nu I\left({x\over {k\over M}\nu}\right)+\eps.
$$
Therefore 
$$
\limsup_{n\ra\infty}{1\over n}\log E[e^{-\lambda H(nx)},\,H(nx)\leq\nu n]\leq-\min_{0<k\leq M}\left({{k\over M}\nu\lambda+{k\over M}\nu I\left({x\over {k\over M}\nu}\right)}\right)+\eps. 
$$
Since the last inequality is true for any $\eps>0$, we obtain 
$$
\limsup_{n\ra\infty}{1\over n}\log E[e^{-\lambda H(nx)},\,H(nx)\leq\nu n]\leq-\inf_{0<\gamma\leq\nu}\left(\gamma\lambda+\gamma I(x/\gamma)\right).
$$
Next, we write 
$$
E[e^{-\lambda H(nx)}]=E[e^{-\lambda H(nx)},\, H(nx)\leq\nu n]\left[1-{E[e^{-\lambda H(nx)},\, H(nx)>\nu n]\over E[e^{-\lambda H(nx)}]}\right]^{-1}.
$$
Using the established lower bound we have that for any $\nu>0$
$$
\limsup_{n\ra\infty}{1\over n}\log{E[e^{-\lambda H(nx)},\, H(nx)>\nu n]\over E[e^{-\lambda H(nx)}]}\leq-\lambda\nu+\inf_{\gamma>0}\left(\gamma\lambda+\gamma I(x/\gamma)\right). 
$$
The last infimum is finite and attained at a finite value $\gamma_0$. Therefore, choosing a value of $\nu$ which satisfies 
$$
\nu>{\gamma_0\over\lambda}(\lambda+I(x/\gamma_0)),
$$
we obtain that 
$$
\lim_{n\ra\infty}{E[e^{-\lambda H(nx)},\, H(nx)>\nu n]\over E[e^{-\lambda H(nx)}]}=0.
$$
Hence 
$$\ea{
\limsup_{n\ra\infty}{1\over n}\log E[e^{-\lambda H(nx)}]=&\limsup_{n\ra\infty}{1\over n}\log E[e^{-\lambda H(nx)},\,H(nx)\leq\nu n]\cr
\leq&-\inf_{0<\gamma\leq\nu}\left(\gamma\lambda+\gamma I(x/\gamma)\right).
}$$
This completes the proof.\qed 
\enddemo
\noi We end this section with a random walk description of the 2d-Ising Wulff crystal.
\proclaim{Corollary \procnum}For all $x=(x_1, x_2)\in\R^2$ and for all $\beta>\beta_c$, the surface tension is given by 
$$
\T(x)=x_1\ac(sx_1)+x_2\ac(sx_2),
$$
where $s$ solves the following equation
$$
\sqrt{1+s^2x_1^2}+\sqrt{1+s^2x_2^2}=\sinh(2\beta)+{1\over\sinh(2\beta)}.
$$
The boundary of the Wulff crystal is given by the following set
$$
\partial\Cal W_\beta=\left\{x\in\R^2:\, \Cal L(x)={1\over 2}\left(\sinh(2\beta)+{1\over\sinh(2\beta)}\right)\right\},
$$
where $\Cal L(x)=E[\exp(x\cdot S_1)]={1\over 2}(\cosh(x_1)+\cosh(x_2))$ is the Laplace transform of the random walk evaluated at $x$.
\endproclaim
\demo{Proof} From Proposition \procref{\tauCramer}, we know that 
$$
\T(x)=\inf_{\gamma>0}\left(\lambda(\beta)\gamma+\gamma I(x/\gamma)\right).
$$ 
Since the function $I$ is even, we have
$$
\T(x)=-\sup_{\gamma<0}(\gamma\lambda(\beta)+\gamma I(x/\gamma)).
$$
The function  $I$ is the Legendre transform of $\Lambda(x)=\log(\Cal L(x))$, so that
$$\ea{
\T(x)=&-\sup_{\gamma<0}\sup_{y\in\R^2}\left((x\cdot y)+\gamma(\lambda(\beta)-\Lambda(y))\right)\cr
=&\sup_{(\gamma, y)}\Phi(y,\gamma),
}$$ 
where 
$$\ea{
\Phi:\,\Omega=(-\infty, 0)\times\R^2&\longrightarrow\R\cr
(\gamma, y)&\longmapsto(x\cdot y)+\gamma(\lambda(\beta)-\Lambda(y)).
}$$
Hence, there exists a point $(\gamma_0, y_0)\in\Omega$ such that $\T(x)=\Phi(\gamma_0, y_0)$ and $(\gamma_0, y_0)$ is the unique solution in $\Omega$ of the equation 
$\nabla\Phi(\gamma, y)=0$. 
Note that the last equation is the Lagrange equation associated to the maximization of $y\mapsto (x\cdot y)$ on the set of $y\in\R^2$ satisfying $\Lambda(y)=\lambda$. Therefore 
$$
\T(x)=\max_{y:\Lambda(y)=\lambda(\beta)}(x\cdot y).\eqnum
$$\eqlabel{\constrvar}%
Since the Wulff shape is the convex polar of the ball corresponding to the norm $\T(x)$, we get that 
$$
\Cal W_\beta=\{x\in\R^2:\q\Lambda(x)\leq\lambda(\beta)\},
$$ 
from which we conclude the statement about the Wulff shape. The formula for the surface tension is obtained by solving the constrained variational problem \eqref{\constrvar}.
\qed
\enddemo
\noi The last corollary describes the Wulff shapes of different temperatures as the level sets of the Laplace transform of the increments of the simple random walk, see the figure below 
\par
\vbox{\line{\hfil\epsfxsize=12cm \epsfbox{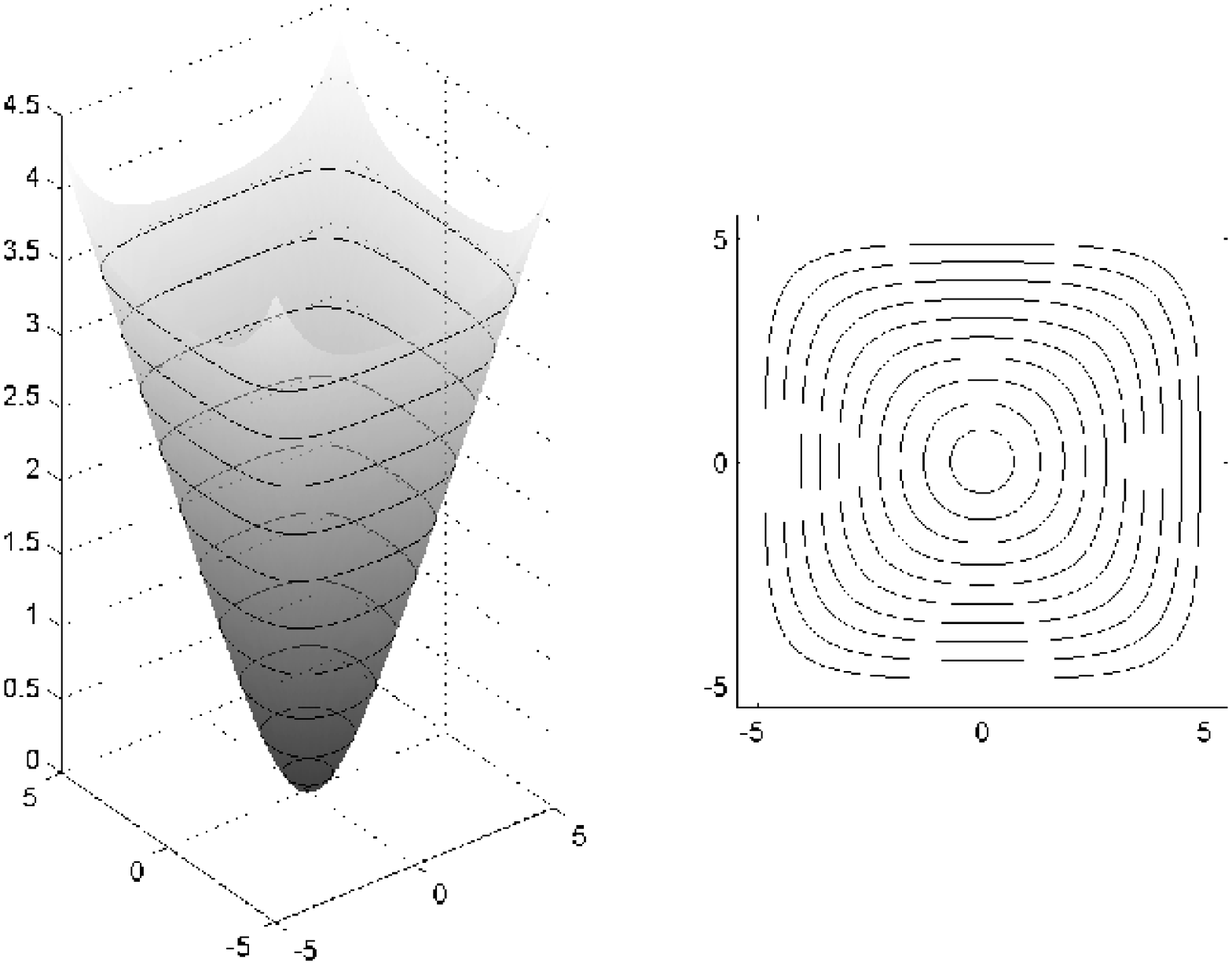}\hfil}\medskip}
\centerline{Left: The graph of $\L(x)\qq$Right: The level sets of $\L(x)$, i.e., the Wulff shapes. }
\medskip\null
\subhead{\headnum The joint limit surface tension}\endsubhead
In this section, we derive the random walk analogue of the surface tension near criticality. 
This will give us a heuristic picture in terms of random walks of Theorem \procref{\jointlimit}. 
It is an interesting question to turn this heuristic description into a rigorous construction of the surface tension in a more general setting than the 2d-Ising model.  

\noi The limit $\beta\da\beta_c$ corresponds, in the random walk picture, to send the surviving probability $m$ to $1$.
Thus, in order to get the analogue of the surface tension near criticality, we consider the asymptotics of $E(\exp(-\lambda H(nx)))$, in the situation where
$n\ua\infty$ and $~{\lambda=-\log(m)\da 0}$ simultaneously. As one may already guess, these asymptotics are related to the moderate deviations of the random walk.  
\proclaim{Proposition \procnum}If $n\ua\infty$ and $\lambda\da 0$
in such a way that
$$
\lim_{n\ua\infty}\frac{1}{n^2\lambda}\log\frac{1}{\lambda}=0,
$$
then
$$ 
\lim_{n\ua\infty,\lambda\da
0}\frac{1}{n\sqrt{\lambda}}\log E[\exp(-\lambda H(nx))]=-2|x|.
$$
\endproclaim\proclabel{\rwjointlimit}
\demo{Proof} From moderate deviations results \cite{\DZ} for the
random walk we can guess the following lower bound:
$$\eqalign{
E[\exp(-\lambda H(nx))]&\geq
E[\exp(-\lambda H(nx)),H(nx)\leq\gamma n/\sqrt{\lambda}]\cr
&\geq\exp(-n\sqrt{\lambda}\gamma)P(S_{\lfloor\gamma n/\sqrt{\lambda}\rfloor}=nx)
}$$   
where $\gamma$ is an arbitrary  positive constant.
We define 
$$
N_n=\frac{\gamma n}{\sqrt{\lambda}}\q,\q a_n=\frac{1}{\gamma n\sqrt{\lambda}}.
$$
Then 
$$
E[\exp(-\lambda H(nx))]\geq\exp(-a_n^{-1})P(S_{\lfloor N_n\rfloor}=\sqrt{\frac{N_n}{a_n}}\
\frac{x}{\gamma}),
$$
so that 
$$
\liminf_{n\ua\infty}\frac{1}{n\sqrt{\lambda}}\log
E(\exp(-\lambda H(nx)))\geq
-\gamma+\gamma\liminf_{n\ua\infty}a_n\log
P(S_{\lfloor N_n\rfloor}=\sqrt{\frac{N_n}{a_n}}\ \frac{x}{\gamma}),$$
Then, from moderate deviations for the simple random walk, as soon as 
$a_n\da 0$ which is equivalent to $n^2\lambda\ua\infty$, 
we obtain 
$$
\liminf_{n\ua\infty}\frac{1}{n\sqrt{\lambda}}\log
E(\exp(-\lambda H(nx)))\geq-\left(\gamma+\gamma\widetilde{I}(x/\gamma)\right),
$$
where $\widetilde{I}(x/\gamma)=|x|^2/\gamma^2$ is the Cram\'er function of the Gaussian
approximation of our random walk.
The upper bound is proved as in Proposition \procref{\tauCramer}. We get 
$$
\liminf_{n\ua\infty}\frac{1}{n\sqrt{\lambda}}\log E(\exp(-\lambda H(nx)))\leq-\inf_{\nu>0}\left(\nu+\nu\widetilde{I}(x/\nu)-\limsup_{n\ua\infty}\frac{-\log\lambda}{n\sqrt{\lambda}}\right).
$$
By imposing $\lim_{n\ua\infty,\lambda\da 0}\log\lambda/n\sqrt{\lambda}=0,$
the lower bound matches the upper bound. The result is obtained by computing the infimum:
$
\inf_{\nu>0}\left(\nu+\nu\widetilde{I}(x/\nu)\right)=2|x|.
$
\qed
\enddemo
\noi Let us stress out that Proposition \procref{\rwjointlimit} does not represent an alternative derivation of Theorem \procref{\jointlimit}. 
It would have been so if we proved that in the joint limit \eqref{\regime} the asymptotic relation \eqref{\formuleMcCoyWuI} is still valid.
Actually, such a proof is possible by adding further restrictions to the regime \eqref{\regime} and it would require a non-trivial modification of the proof of Mc Coy and Wu \cite{\McCoyWu}. This would add another layer of explicit computations. 
\subhead\subheadnum The explicit computations formula and the method of Laplace\endsubhead
Until now we considered only rough asymptotics of $V_m(x_1, x_2)$ that were good enough to describe the surface tension. 
In this section, we use the method of Laplace on the double integral to reproduce the prefactor in the results of Mc Coy and Wu. 
This will allow us to view the Ornstein-Zernike behavior as a result of a Laplace method. 
First, we express $V_m(x_1,x_2)$ in terms of Bessel functions. 
For an extensive treatment of these special functions, we refer the reader to \cite{\WATSON}.
\proclaim{Lemma \procnum} For every $(x_1,x_2)\in\Z^2$ we have that 
$$\eqalign{V_m(x_1,x_2)=&
\frac{r}{\pi}\frac{x_1^{x_1}x_2^{x_2}}{\Gamma(x_1+1/2)\Gamma(x_2+1/2)}\cr
&\times\int_{[0,+\infty)}du\int_{[0,\pi]}d\theta_1\int_{[0,\pi]}d\theta_2\exp(rf_m(u,\theta_1,\theta_2)),
}\eqnum$$\eqlabel{\massdeux}
\def\expfm{&\frac{mu}{2}(\cos\theta_1+\cos\theta_2)-u\cr
\textstyle&+\cos\varphi\log\left(\frac{m\,u\,\sin^2\theta_1}{4\cos\varphi}\right)+\sin\varphi\log\left(\frac{m\,u\,\sin^2\theta_2}{4\sin\varphi}\right)}where
$$\eqalign{
f_m(u,\theta_1,\theta_2)=\expfm.
}$$
and $r=\sqrt{{x_1}^2+{x_2}^2}$. 
\endproclaim
\proclabel{\formulevisites}
\demo{Proof}
Noticing that each walk that starts at the origin and ends at $(x_1,x_2)$ has to contain $x_1+2k_1$ horizontal displacements and $x_2+2k_2$ vertical displacements with $k_1$ and  $k_2$ are two positive integers, we get that 
\def\LR#1{\left(#1\right)}
{\def\t{x_1+x_2+2(k_1+k_2)}
$$\eqalign{
V_m(x_1,x_2)&=\sum_{k_1,k_2=0}^\infty m^{\t} P(X_{\t}=(x_1,x_2))\cr
&=\sum_{k_1,k_2=0}^\infty \frac{(\t)!}{k_1!k_2!(x_1+k_1)!(x_2+k_2)!}\LR{\frac{m}{4}}^{\t}.
}$$
Using the identity $\int_0^\infty du\,e^{-r u}u^n=n!\;r^{-n-1}$, we obtain
$$\eqalign{
V_m(x_1,x_2)&=r\int_{0}^{\infty}du\;e^{-r u}\sum_{k_1=0}^\infty\frac{\LR{m r u/4}^{x_1+2k_1}}{k_1!(x_1+k_1)!}\sum_{k_2=0}^\infty\frac{\LR{m r u/4}^{x_2+2k_2}}{k_2!(x_2+k_2)!}
.}$$
}%
Next, we recognize the expansion of the Bessel function of the first type \cite{\WATSON}
$$
\forall n\in\N\quad\forall x\in\R\qquad\,I_n(x )=\sum_{k=0}^{\infty}\frac{1}{k!(n+k)!}\LR{\frac{x }{2}}^{n+2k}.
$$
and get
$$V_m(x_1,x_2)=r\int_0^\infty
e^{-ru}I_{x_1}(\frac{m\,r}{2}u)I_{x_2}(\frac{m\,r}{2}u)\,du.\eqnum
$$\eqlabel{\massun}%
Using the following representation of Bessel functions \cite{\WATSON}
\medskip
$$
I_n(\alpha n)=\frac{\LR{\alpha
n/2}^n}{\sqrt{\pi}\;\Gamma(\frac{1}{2}+n)}\int_0^\pi d\theta\,
\exp(n(\alpha\cos\theta+\log\sin^2\theta)),
$$
we write
$$\eqalign{
I_{x_1}\Bigl(\frac{mru}{2}\Bigl)&=\frac{\Bigl(\frac{mux_1}{4\cos\vphi}\Bigr)^{x_1}}{\sqrt{\pi}\;\Gamma(x_1+\frac{1}{2})}\int_0^\pi
d\theta
\exp\left[r\LR{\frac{mu}{2}\cos\theta+\cos\vphi\log\sin^2\theta}\right]\cr
&=\frac{x_1^{x_1}}{\sqrt{\pi}\;\Gamma(x_1+\frac{1}{2})}\int_0^\pi
d\theta
\exp\left[r\LR{\frac{mu}{2}\cos\theta+\cos\vphi\log\frac{mu\sin^2\theta}{4\cos\vphi}}\right].
}$$ 
and in a similar way
$$
I_{x_2}\Bigl(\frac{mru}{2}\Bigl)=\frac{x_2^{x_2}}{\sqrt{\pi}\;\Gamma(x_2+1/2)}\int_0^\pi
d\theta\exp\left[r\LR{\frac{mu}{2}\cos\theta+\sin\vphi\log\frac{mu\sin^2\theta}{4\sin\vphi}}\right].
$$
The result follows by replacing the two last expressions for the Bessel functions into \eqref{\massun}. 
\qed
\enddemo
\noi Lemma \procref{\formulevisites} gives a convenient way to treat the asymptotics with the method of Laplace. 
We will use the following version of this method and refer the reader to \cite{\HSU}, for its derivation.  

\proclaim{Proposition \procnum} Let  $D\subset\R^d$ be an open set and consider a function
$f:D\longrightarrow\R$ that satisfies
\medskip
{\parindent=10mm
\item{i)} {$\forall r>0\quad\exp(r f(x))$ is integrable on $D$.}
\item{ii)} {$f$ is twice differentiable on $D$.}
\item{iii)} {$f$ reaches its global maximum on $D$.}
\item{iv)} {The maximum of $f$ is reached at a unique point ${\xs}\in D\setminus\partial D$.}
\item{v)}{$\forall x_0\in\partial D\quad\lim_{x\ra x_0}f(x)=-\infty.$}
\par}
\medskip\noindent Then one has
$$ \int_D\exp(r
f(x))\,dx\sim\left(\frac{2\pi}{r}\right)^{d/2}\frac{\exp(r
f(\xs))}{\sqrt{\det(-\H(\xs))}}\quad r\ra\infty, \eqnum$$
\eqlabel{\asympLaplace}where $\H(\xs)$ is the Hessian of $f$ at the critical point $\xs$.
\endproclaim
\proclabel{\methLaplace}

\noi Now we are ready to compute the asymptotics of \eqref{\massdeux}, in a situation where $m$ is fixed and $r\ra\infty$.
\proclaim{Proposition \procnum}
Let $\beta<\beta_c$ be fixed. For every $\varphi\in[0,\pi/4]$, the asymptotics of $\mu^{\beta}[\sigma(0)\sigma(r\cos\varphi, r\sin\varphi)]$ when $r\ra\infty$ are given by 
$$\ea{
&
{(\sinh^{-4}(2\beta)-1)^{1/4}\over(1-\tanh^4(\beta))^2\sinh(2\beta)}
\frac{\sqrt{{\us/\pi mr}}}{\left(\cos^2\vphi\sqrt{\frac{4\sin^2\vphi}{m^2\us^2}+1}+\sin^2\vphi\sqrt{\frac{4\cos^2\vphi}{m^2\us^2}+1}\right)^{1/2}}\times\cr
&\times\exp\LR{-r\LR{\cos\vphi\;\ac\frac{2\cos\vphi}{m\us}+\sin\vphi\;\ac\frac{2\sin\vphi}{m\us}}},
}\eqnum$$
\eqlabel{\AsympVisites}where
$$
\us= \left(\frac{\left((1-m^2)\sin^22\vphi+m^2\right)^{1/2}+1}{1-m^2}\right)^{1/2},
$$
and
$$
m={2\over\sinh(2\beta)+{1\over\sinh(2\beta)}}.
$$
\endproclaim
\proclabel{\asympr}
\demo{Proof}
We apply the method of Laplace given by Proposition \procref{\methLaplace} on the function
$$\eqalign{
\textstyle f_m(u,\theta_1,\theta_2)=\expfm,
}$$
defined on $D=(0,+\infty)\times(0,\pi)^2$.

\noi To show that $f_m$ reaches its global maximum at a unique point, we note that for every $(\theta_1,\theta_2)\in[0,\pi]^2$
one has
 $$-(1-\frac{m}{2}(\cos\theta_1+\cos\theta_2))\leq-(1-m)<0.$$
 Thus for $x_0\in\partial
D$ we certainly have  that $\lim_{x\ra x_0}f_m(x)=-\infty$ because
$$\lim_{(\theta_1,\theta_2)\rightarrow(0,0)}f_m(u,\theta_1,\theta_2)=-\infty,
$$uniformly  in $u$ and  
$$\lim_{u\rightarrow
+\infty}f_m(u,\theta_1,\theta_2)=-\infty,$$ uniformly in $\theta_1$ et $\theta_2$. Hence $f_m$ reaches its maximum in  $D\setminus\partial D$ at a critical point that satisfies

$$
\grad f_m(u,\theta_1,\theta_2)=0\Leftrightarrow
\left\{
\matrix
-1+\frac{m}{2}(\cos\theta_1+\cos\theta_2)+\frac{\cos\vphi}{u}+\frac{\sin\vphi}{u} &=& 0\cr
-\frac{m}{2} u\sin\theta_1+2\cos\vphi\cot\theta_1&=&0\cr
-\frac{m}{2} u\sin\theta_2+2\sin\vphi\cot\theta_2&=&0\cr
\endmatrix
\right. $$ 
It is easy to solve these equations and to get that it admits a unique solution $(\us, \ths_1, \ths_2)$ given by
$$ \left\{\eqalign{
\us=\us_+=&\sqrt{\frac{1+\sqrt{1-(1-m^2)\cos^22\vphi}}{1-m^2}}\cr
\cos\ths_1=&\Exp{-\ac\left(\frac{2\cos\vphi}{m\,\us}\right)}=\sqrt{\LR{\frac{2\cos\vphi}{m\,\us}}^2+1}-\frac{2\cos\vphi}{m\,\us}\cr
\cos\ths_2=&\Exp{-\ac\left(\frac{2\sin\vphi}{m\,\us}\right)}=\sqrt{\LR{\frac{2\sin\vphi}{m\,\us}}^2+1}-\frac{2\sin\vphi}{m\,\us}\cr
} \right.$$
The Hessian matrix $\H(u, \theta_1, \theta_2)$ of $f_m$ at a point $(u,\theta_1, \theta_2)$ is given by 
{\hfuzz=20pt
$$
\pmatrix
-\frac{\cos\vphi+\sin\vphi}{u^2} &  -\frac{m}{2}\sin\theta_1  &  -\frac{m}{2}\sin\theta_2\cr
\cr
-\frac{m}{2}\sin\theta_1  &  -\frac{m}{2}u\cos\theta_1-2\frac{\cos\vphi}{\sin^2\theta_1}  &  0\cr
\cr
-\frac{m}{2}\sin\theta_2  &  0  &  -\frac{m}{2} u\cos\theta_2+2\frac{\sin\vphi}{\sin^2\theta_2}
\endpmatrix
$$}

\noi and at the critical point, the determinant of $\Hs=\H(\us, \ths_1, \ths_2)$ can be computed, we get
$$
\det\Hs=-\frac{2m}{\us}\left(\cos^2\vphi\sqrt{\left(\frac{2\sin\vphi}{m\us}\right)^2+1}+\sin^2\vphi\sqrt{\left(\frac{2\cos\vphi}{m\us}\right)^2+1}\right)<0.
\eqnum$$
\eqlabel{\ValDetExt}Thus, the method of Laplace can be applied. For this, we compute the maximum of  $f_m$  and get
$$
f_m(\us,\ths_1,\ths_2)=-\cos\vphi\,\ac\frac{2\cos\vphi}{m\us}-\sin\vphi\,\ac\frac{2\sin\vphi}{m\us}-\cos\vphi-\sin\vphi\eqnum
$$
\eqlabel{\ValFoncExt}The asymptotics of the integral in $(\massdeux)$ are then obtained by replacing \eqref{\ValDetExt} and \eqref{\ValFoncExt} into \eqref{\asympLaplace}.
We get that for each $0<m<1$ 

$$\eqalign{
&\int_{[0,+\infty]\times[0,\pi]^2}\exp(rf_m(x))\,dx\sim\left(\frac{2\pi}{r}\right)^{3/2}\frac{\exp(r f_m(\xs))}{\sqrt{\det(-\Hs)}}\cr
&\cr
=&\frac{\exp\LR{-r\LR{\cos\vphi\;\ac\frac{2\cos\vphi}{m\us}+\sin\vphi\;\ac\frac{2\sin\vphi}{m\us}}}e^{r(\cos\vphi+\sin\vphi)}}{\left(\frac{r}{2\pi}\right)^{3/2}\LR{\frac{2m}{\us}}^{1/2}\left(\cos^2\vphi\sqrt{\frac{4\sin^2\vphi}{m^2\us^2}+1}+\sin^2\vphi\sqrt{\frac{4\cos^2\vphi}{m^2\us^2}+1}\right)^{1/2}},
}\eqnum$$
when $r\ra\infty$.
\eqlabel{\asympint}
Finally, we treat the prefactor 
$x_1^{x_1}x_2^{x_2}/\pi\Gamma(x_1+1/2)\Gamma(x_2+1/2)$. From Stirling's formula, we get 
$$\frac{x_1^{x_1}}{\Gamma(x_1+1/2)}\frac{x_2^{x_2}}{\Gamma(x_2+1/2)}\frac{r}{\pi}\sim\frac{r e^{r(\cos\vphi+\sin\vphi)}}{2\pi^2},\quad r=\norm{(x_1,x_2)}\ra\infty.
\eqnum
$$\eqlabel{\asympfacteur}

\noi Combining \eqref{\asympfacteur}, \eqref{\asympint} and \eqref{\formuleMcCoyWuII}
, we get the desired result.
\qed
\enddemo

\enddocument